\documentclass[12pt]{article}

\usepackage{amsmath}
\usepackage{amssymb}
\usepackage{amsfonts}
\usepackage{graphics}

\newtheorem{thm}{Theorem}[section]

\newtheorem{prop}[thm]{Proposition}
\newtheorem{conjecture}{Conjecture}

\newcommand{\be}{\begin{equation}}
\newcommand{\ee}{\end{equation}}
\newcommand{\openbox}{\leavevmode
  \hbox to8pt{\hfil\vrule\vbox to6pt{\hrule width6pt\vfil\hrule}\vrule}}

\newcommand{\qed}{\hbox to5pt{ } \hfill \openbox\bigskip\medskip}

\newcommand{\cF}{\mbox{$\cal F$}}

\newcommand{\ve}[1]{\mathbf{#1}}

\newcommand{\R}{\mathbb R}

\newcommand{\F}{\mathbb F}

\newcommand{\rank}{\mathop\textup{rank}}

\title{A uniform version of a theorem by Lindstr\"om }
\author{G\'abor Heged\"{u}s
\\{\normalsize  \'Obuda University}
\\{\normalsize B\'ecsi \'ut 96/B, Budapest, Hungary, H-1032}
\\{\normalsize hegedus.gabor@uni-obuda.hu}
}

\begin{document}

\maketitle
\begin{abstract}
We prove the following uniform version of a theorem by Lindstr\"om: \\
Let  $\mbox{$\cal F$}:=\{F_i:~ i\in I\}$ be a $k$-uniform set family on $[n]$, where $k\geq 1$. If $|\mbox{$\cal F$}|\geq n+1$, then there exist two disjoint subsets $I_1$ and $I_2$ of $I$ for which
$$
\bigcup\limits_{i\in I_1} M_i=\bigcup\limits_{i\in I_2} M_i 
$$
and 
$$
\bigcap\limits_{i\in I_1} M_i=\bigcap\limits_{i\in I_2} M_i.
$$

Our proof uses basic linear algebra.
\end{abstract}
\medskip

\footnotetext{ 

 05D05; 15A06; 15A03}

\medskip

\section{Introduction}

Let $n$ denote a positive integer. Let $[n]$ stands for 
the set $\{1,2,\dots, n\}$ and  denote 
by $2^{[n]}$ the set of all subsets of $[n]$.
The subsets of $2^{[n]}$ are called simply \emph{set families}.
We   denote  by $\binom{[n]}m$ the family of all subsets of $[n]$
which have cardinality $m$.

A family $\cF$ of subsets of $[n]$ is {\em $k$-uniform}, if $|F|=k$ for each $F\in \cF$.


Consider a family $\cF=\{F_i:~ i\in I\}$ of subsets of $[n]$. We say that $\cF$ is a {\em $\cup$-balanced  family}, if we find two disjoint subsets $I_1$ and $I_2$ of $I$ for which
$$
\bigcup\limits_{i\in I_1} M_i=\bigcup\limits_{i\in I_2} M_i.
$$

A family  $\cF=\{F_i:~ i\in I\}$ of subsets of $[n]$ is a  {\em balanced  family}, if there exist two disjoint subsets $I_1$ and $I_2$ of $I$ for which
$$
\bigcup\limits_{i\in I_1} M_i=\bigcup\limits_{i\in I_2} M_i 
$$
and 
$$
\bigcap\limits_{i\in I_1} M_i=\bigcap\limits_{i\in I_2} M_i.
$$

To the best of the author's knowledge,  the following result is folklore.
\begin{thm}\label{lind1}
Let  $\cF$ be a set family on $[n]$. If $|\cF|\geq n+1$, then $\cF$ is a $\cup$-balanced  family.
\end{thm}

Lindstr\"om gave a generalization of Theorem  \ref{lind1} in \cite{L1} using Rado's extension of Hall's Theorem. In \cite{L2} Lindstr\"om proved the following result for balanced set families using linear algebra and asked for a purely combinatorial proof.  

\begin{thm}\label{lind2}
Let  $\cF$ be a set family on $[n]$. If $|\cF|\geq n+2$, then $\cF$ is a balanced  family.
\end{thm}
It is easy to verify that the above inequality is sharp.
Namely let $\cF:=\{\{i\}:~ 1\leq i\leq n\}\cup \{[n]\}$. Then $|\cF|= n+1$ and $\cF$ is not  a balanced family.
Motivated by this result we investigated uniform balancing families and we prove here the following uniform version of Theorem \ref{lind2}. 

\begin{thm}\label{main}
Let  $\cF:=\{F_1, \ldots ,F_m\}$ be a $k$-uniform set family on $[n]$, where $k\geq 1$. If $|\cF|\geq n+1$, then $\cF$ is a balanced  family.
\end{thm}
The following set system shows that the above inequality is sharp: let $\cF:=\{\{2,3\}\}\cup \{\{1,i\}:~ 2\leq i\leq n\}$. Clearly $|\cF|= n$, $\cF$ is uniform and $\cF$ is not balanced.

Our proof  uses the following simple Proposition (see \cite{BF} Proposition 2.12).
\begin{prop} \label{rank}
Let $\mathbf{S}$ be an $m\times n$ matrix over a field $\F$. The set of solutions of the system of linear equations $\mathbf{S}\cdot \ve x=0$ is a linear subspace of dimension $n- \rank(\mathbf{S})$ of the space $\F^n$.
\end{prop}

\section{The proof of our main result}
 
Let $V\subseteq \R^{2n}$ denote the real vector space of all vectors  $\ve v=(x_1, y_1, \ldots x_n, y_n)$ for which $x_i+y_i=x_j+y_j$ for each $1\leq i\neq j\leq n$ and $(n-k)\sum_{i=1}^n x_i=k\sum_{i=1}^n  y_i$.

We claim that $\dim V=n$.

Namely the space $V$ is precisely the set of solutions $\ve v\in \R^{2n}$ of the system of linear equations $\mathbf{T}\cdot \ve v=0$, where $T$  is the following $n\times 2n$ matrix:
\[ \mathbf{T}=\left(
\begin{array}{ccccccccc}
1 & 1 & -1 & -1 &  0 & 0 & \cdots & 0 & 0 \\
1 & 1 & 0 & 0 &  -1 & -1 & \cdots & 0 & 0 \\
\vdots & \vdots & \vdots & \vdots & \vdots & \vdots & \ddots & \vdots & \vdots \\
1 & 1 & 0 & 0 &  0 & 0 & \cdots &  -1 & -1 \\
(n-k) & -k & (n-k) & -k & (n-k) & -k & \cdots &  (n-k) & -k 
\end{array} \right) \]

It is easy to verify that $\rank(\mathbf{T})=n$, hence $\dim V=2n-\rank(\mathbf{T}) =n$ by Proposition \ref{rank}. 

Let $\overline{A}$ denote the complement of the subset $A$. We associate with each subset $A$  the {\em extended incidence vector} $\ve v=(x_1, y_1, \ldots x_n, y_n)$ of the pair $(A,\overline{A})$ the following way: let  $x_i=1$ iff $i\in A$ (otherwise $x_i=0$)  and define $y_i:=1-x_i$ for each $1\leq i \leq n$. 

Let  $\ve v_i$ denote the extended incidence vector, which corresponds  to the $i$th set $A_i$,  where $1\leq i \leq m$. It is easy to verify from the assumptions of Theorem \ref{main} that the vectors $\ve v_i$ are distinct and $\ve v_i\in V$ for each $1\leq i \leq m$. It follows from  $m>\dim V$, that  there must be a nontrivial  linear relation between the vectors  $\ve v_i\in V$. Clearly we can  write this relation as
\begin{equation} \label{sum}
\sum_{i\in I_1} \alpha_i\ve v_i =\sum_{j\in I_2} \beta_j\ve v_j,
\end{equation}
where $I_1$ and $I_2$ are non-empty, disjoint sets, and $ \alpha_i,\beta_j>0$ for all $i\in I_1, j\in I_2$. The equality (\ref{sum}) implies  that
$$
\bigcup\limits_{i\in I_1} A_i=\bigcup\limits_{j\in I_2} A_j
$$
and
$$
\bigcup\limits_{i\in I_1} \overline{A_i}=\bigcup\limits_{j\in I_2} \overline{A_j}.
$$
It follows from de Morgan's law that the last equality is equivalent to 
$$
\bigcap\limits_{i\in I_1} A_i=\bigcap\limits_{j\in I_2} A_j.
$$ 
\qed

\section{Concluding remark}

We propose the following conjecture as a natural strengthening of
Theorem \ref{main}.
\begin{conjecture} \label{Hconj}
Let  $\cF=\{F_1, \ldots ,F_m\}$ be a Sperner set family on $[n]$. If $|\cF|\geq n+1$, then $\cF$ is a balanced  family.
\end{conjecture}

\end{document}